\documentclass{amsart}

\def\ra{\rightarrow}
\def\bs{\hskip -.1in}
\def\Sb#1{_{\substack{#1}}}

\def\Z{\mbox{$\mathbb Z$}}
\def\N{\mbox{$\mathbb N$}}

\def\ph{\varphi}

\def\and{\hbox{ and }}
\def\s#1{\mbox{$(-1)^{#1}$}}
\def\e#1{\mbox{$|#1|$}}
\def\ip#1#2{\left<#1,#2\right>}
\def\br#1#2{\left\{#1,#2\right\}}
\def\brf{\left\{\cdot,\cdot\right\}}
\def\brt#1#2{\left\{#1,#2\right\}}
\def\brtf{\left\{\cdot,\cdot\right\}}

\def\ipf{\left<\cdot,\cdot\right>}
\def\blacksquare{\quad \vrule height 6pt width 6pt}
\def\inv{^{-1}}

\newtheorem{thm}{Theorem}

\newtheorem{lma}{Lemma}
\def\ie{\hbox{\it i.e.}}
\def\ainf{\hbox{$A_\infty$}}
\def\linf{\mbox{$L_\infty$}}
\def\zt{\hbox{$\Z_2$}}
\def\ztz{\hbox{$\Z_2\times\Z$}}
\def\hd{\hbox{{$\hat d$}}}
\def\hm{\hbox{$\hat m$}}
\def\hl{\hbox{$\hat l$}}
\def\hmu{\hbox{$\hat \mu$}}

\def\tf{\hbox{$\tilde f$}}

\def\bmu{\hbox{$\bar\mu$}}

\def\bm{\hbox{$\bar m$}}

\def\hdl{\hbox{$\hat \delta$}}

\def\CW{C(W)}
\def\CV{C(V)}
\def\hom{\mbox{\rm Hom}}

\def\sb#1{[#1]}%superbracket notation
\def\tens{\bigotimes}

\def\tns{\otimes}
\def\mtns{\tns\cdots\tns}
\def\modot{\odot\cdots\odot}

\def\mwedge{\wedge\cdots\wedge}
\def\modot{\odot\cdots\odot}
\def\mplus{+\cdots+}
\def\mcom{,\cdots,}
\def\iso{\kern.35em{\raise3pt\hbox{$\sim$}\kern-1.1em\to}
         \kern.3em}

\def\k{\mbox{\bf k}}

\def\V{\mbox{V}}
\def\SV{S(V)}
\def\EV{\bigwedge(V)}
\def\TV{T(V)}

\def\bid{\operatorname{bid}}
\def\bd#1{\bid(#1)}
\def\shf{\operatorname{Sh}}
\def\sh#1#2{\shf(#1,#2)}

\def\deg{\operatorname{deg}}

\def\coder{\operatorname{Coder}}

\def\s#1{\mbox{$(-1)^{#1}$}}
\def\e#1{\mbox{$|#1|$}}
\def\ip#1#2{\left<#1,#2\right>}
\def\blacksquare{\quad \vrule height 6pt width 6pt}

\begin{document}
\nocite{mar,sta3,getz2,lod,conn,kast,seib,hoch,umb,aksz}
\author{Michael Penkava}
\thanks{Partially supported by NSF grant DMS-94-0411.
The author would also like to thank the University of Washington,
Seattle for hosting him.}
\address{University of California\\
Davis, CA 95616}
\email{michae@@math.ucdavis.edu}
\title[Cohomology of Infinity Algebras]{Infinity Algebras, Cohomology and Cyclic Cohomology,
and Infinitesimal Deformations}
\begin{abstract}
% Abstract text begins here
An $A_\infty$ algebra is given by a codifferential on the tensor
coalgebra of a (graded) vector space. An associative algebra is a
special case of an $A_\infty$ algebra, determined by a quadratic
co\-differential. The notions of Hochschild and cyclic cohomology
generalize from associative to $A_\infty$ algebras, and classify
the infinitesimal deformations of the algebra, and those
deformations preserving an invariant inner product, respectively.
Similarly, an $L_\infty$ algebra is given by a codifferential on
the exterior coalgebra of a vector space, with Lie algebras being
special cases given by quadratic codifferentials. There are
natural definitions of cohomology and cyclic cohomology,
generalizing the usual Lie algebra cohomology and cyclic
cohomology, which classify deformations of the algebra and those
which preserve an invariant inner product.  This article explores
the definitions of these infinity algebras, their cohomology and
cyclic cohomology, and the relation to their infinitesimal
deformations.
\end{abstract}
\maketitle
\section{Introduction}

In a joint paper with Albert Schwarz \cite{ps2}, we gave
definitions of Hochschild cohomology and cyclic cohomology of an
\ainf\ algebra, and showed that these cohomology theories
classified the infinitesimal deformations of the \ainf\ structure
and those deformations preserving an invariant inner product,
respectively. Then we showed that the Hochschild cohomology of an
associative algebra classifies the deformations of the algebra
into an \ainf\ algebra, and the cyclic cohomology of an algebra
with an invariant inner product classifies the deformations of the
algebra into an \ainf\ algebra preserving  the inner product.
Following some ideas of Maxim Kontsevich from \cite{kon}, we
applied these results to show that cyclic cocycles of an
associative algebra determine homology cycles in the complex of
metric ribbon graphs. In \cite{ps2}, for sake of simplicity we
avoided the more standard description of \ainf\ algebras in terms
of codifferentials.

In this article, the definitions of infinity algebras will be
given in terms of co\-differentials on coalgebras, so that the
theory of \linf\ algebras is seen to be closely analogous to that
of \ainf\ algebras. The ideas in this text lead immediately to a
simple formulation of the cycle in the complex of metric ribbon
graphs  associated to an \ainf\ algebra with an invariant inner
product. This same method can be applied to show that \linf\
algebras with an invariant inner product give rise to a cycle in
the homology of the complex of metric ordinary graphs (see
\cite{pen3}). This paper is a revision of \cite{pen2}.

The notion of an \ainf\ algebra, also called a strongly homotopy
associative algebra, was introduced by J. Stasheff in
\cite{sta1,sta2}, and is a generalization of an associative
algebra.  From a certain point of view, an associative algebra is
simply a special case of a codifferential on the tensor coalgebra
of a vector space. An \ainf\ algebra is given by taking an
arbitrary coderivation; in particular associative algebras  and
differential graded associative algebras are examples of \ainf\
algebras.

\linf\ algebras, also called strongly homotopy Lie algebras, first
appeared in \cite{ss}, and are generalizations of Lie algebras. A
Lie algebra can be viewed as simply a special case of a
codifferential on the exterior coalgebra of a vector space, and
\linf\ algebras are simply arbitrary codifferentials on this
coalgebra.

A bracket structure was introduced on the space of cochains of an
associative algebra by Murray Gerstenhaber in \cite{gers}. {}From
the coalgebra point of view, the \emph{Gerstenhaber bracket} turns
out to be essentially the commutator of coderivations, as was
first pointed out by James Stasheff in \cite{sta4}. In the
presence of an invariant inner product, the cyclic cochains can be
identified as a subalgebra of the space of ordinary cochains, so
that the cyclic cohomology also is equipped with a natural Lie
bracket. These notions generalize to the \ainf\ case.

The space of cochains of a Lie algebra with coefficients in the
adjoint representation has a natural bracket, which is again the
commutator of coderivations. In the case of a Lie algebra, cyclic
cohomology corresponds to the cohomology of the Lie algebra with
trivial coefficients, and in the presence of an invariant inner
product, the space of cyclic cochains is also equipped with a
natural bracket. These notions generalize to the \linf\ algebra
case.

In our considerations, we shall be interested in \zt-graded
spaces, but we should point out that the results presented here
hold in the \Z-graded case as well, because the signs in the
\Z-graded case coincide with the signs from the associated
\zt-grading. \ainf\ algebras were first defined as \Z-graded
objects, but for the applications we have in mind, the \zt-grading
is more appropriate, and the generalization of the results here to
the \Z-graded case is straightforward. We shall find it necessary
to consider the parity reversion of a \zt-graded space. This is
the same space with the parity of elements reversed. (In the
\Z-graded case, the corresponding notion is that of suspension.)

There is a natural vector space isomorphism between the tensor
coalgebra of a \zt-graded space and the tensor coalgebra of its
parity reversion. But in the case of the exterior coalgebra, the
isomorphism is to the {\em symmetric} coalgebra of the parity
reversion.

A notion that will play a crucial role in what follows is that of
a grading group. An abelian group $G$ is said to be a grading
group if it possesses a symmetric \zt-valued bilinear form $\ipf$.
Any abelian group with a subgroup of index 2 possesses a natural
grading form. An element $g$ of $G$ is called odd if $\ip gg=1$. A
grading on a group will be called {\em good} if $\ip gh=1$,
whenever $g$ and $h$ are both odd. Groups equipped with the
natural inner product induced by a subgroup of index 2 are good,
and these include both \zt\ and \Z. If $G$ and $H$ are grading
groups, then $G\times H$ has an induced inner product, given by
$\ip{(g,g')}{(h,h')}=\ip g{g'}+\ip h{h'}$. But the induced grading
on $G\times H$ is never good when the gradings on $G$ and $H$ are
good.

If $V$ is a $G$-graded vector space, then one can define the
symmetric and exterior algebras of the tensor algebra $T(V)$. The
symmetric algebra is $G$-graded commutative, but the exterior
algebra is not. On the other hand, the tensor, symmetric and
exterior algebras also are graded by $G\times\Z$, and with respect
to the induced inner product on $G\times\Z$, the exterior algebra
is graded commutative. Thus the natural grading associated to the
exterior algebra is not good. The consequences of this fact
complicates the definition of \ainf\  and \linf\ algebras.
In the following, we shall always consider the case $G=\zt$,
although the case $G=\Z$ is similar.
\section{The Exterior and Symmetric Algebras}\label{sect 3}

Suppose that $V$ is a \zt-graded \k-module. The (reduced) tensor
algebra $T(V)$ is given by $T(V)=\bigoplus_{n=1}^\infty V^n$,
where $V^n$ is the $n$-th tensor power of $V$. The full tensor
algebra is defined to include the term $V^0=\k$ as well, and in
\cite{getz} a notion of an \ainf\ algebra on the full tensor
algebra is given, but we do not treat this idea here.

For an element $v=v_1\mtns v_n$ in $T(V)$, define its parity
$\e{v}=\e{v_1}\mplus \e{v_n}$, and its degree by $\deg(v)$=n. We
define the bidegree of $v$ by $\bd v=(\e{v},\deg(v))$. If $\alpha,
\beta \in T(V)$, then $\bd{\alpha\tns\beta}=\bd \alpha+\bd\beta$,
so that $T(V)$ is naturally $\ztz$-graded by the bidegree, and
\zt-graded if we consider only the parity. For simplicity, for a
\zt-graded space, we will denote
$\s{\alpha\beta}=\s{\e\alpha\e\beta}$, and for a \ztz-graded space
(equipped with the product grading), we shall denote
$\s{\ip\alpha\beta}=\s{\ip{\bid\alpha}{\bid\alpha}}$. (In formulas
like the above, the elements $\alpha$ and $\beta$ are assumed to
be homogeneous.)

The (graded) symmetric algebra $\SV$ is defined as the quotient of
the tensor algebra of $V$ by the bigraded ideal generated by all
elements of the form $u\tns v-\s{uv}v\tns u$. The resulting
algebra has a decomposition $\SV=\bigoplus_{n=1}^\infty S^n(V)$,
where $S^n(V)$ is the image of $V^n$ in $\SV$. We shall denote the
induced product in $\SV$ by juxtaposition. The symmetric algebra
is both \zt\ and $\ztz$-graded, but is graded commutative only
with respect to the \zt-grading.

If $\sigma\in\Sigma_n$, then
\begin{equation}
v_1\cdots v_n=\epsilon(\sigma;v_1\mcom v_n)v_{\sigma(1)}\cdots
v_{\sigma(n)}
\end{equation}
where $\epsilon(\sigma;v_1\mcom v_n)$ is a sign which can be
determined by the following. If $\sigma$ interchanges $k$ and
$k+1$, then $\epsilon(\sigma;v_1\mcom v_n)=\s{v_kv_{k+1}}$. In
addition, if $\tau$ is another permutation, then
\begin{equation}
\epsilon(\tau\sigma;v_1\mcom v_n)=
\epsilon(\sigma;v_{\tau(1)}\mcom v_{\tau(n)})
\epsilon(\tau;v_1\mcom v_n).
\end{equation}
It is conventional to abbreviate $\epsilon(\sigma;v_1\mcom v_n)$
as $\epsilon(\sigma)$.

The (graded) exterior algebra $\bigwedge V$ is defined as the
quotient of $T(V)$ by the bigraded ideal generated by all elements
of the form $u\tns v+\s{uv}v\tns u$. The resulting algebra has a
decomposition $\bigwedge V=\bigoplus_{n=0}^\infty \bigwedge^n V$,
and the induced product is denoted by $\wedge$. The exterior
algebra is both \zt\ and $\ztz$-graded, but is graded commutative
only with respect to the $\ztz$-grading. Let $\alpha,\beta\in
\bigwedge V$. Furthermore it is easy to see that
\begin{equation}
v_1\mwedge v_n=\s{\sigma}\epsilon(\sigma)v_{\sigma(1)}\mwedge
v_{\sigma(n)},
\end{equation}
where $\s\sigma$ is the sign of the permutation $\sigma$.
\section{The Tensor, Exterior and Symmetric Coalgebras}\label{sect
4} A more formal treatment of the symmetric and exterior
coalgebras would introduce the coalgebra structure on the tensor
algebra, and then describe these coalgebras in terms of coideals.
Instead, we will describe these coalgebra structures directly.
Recall that a coalgebra structure on a \k-module $C$ is given by a
diagonal mapping $\Delta:C\ra C\tens C$. We consider only
coassociative coalgebras (we do not consider counits). The axiom
of coassociativity is
\begin{equation*}(1\tns \Delta)\circ
\Delta=(\Delta\tns 1)\circ\Delta.
\end{equation*}
A grading on $C$ is compatible with the coalgebra structure if for
homogeneous $c\in C$, $\Delta(c)=\sum_i u_i\tens v_i$, where
$\e{u_i}+\e{v_i}=\e c$ for all $i$; in other words, $\Delta$ has
degree 0. We also mention that a coalgebra is graded cocommutative
if $S\circ\Delta=\Delta$, where $S:C\tens C\ra C\tens C$ is the
symmetric mapping given by
\begin{equation*}
S(\alpha\tns\beta)=\s{\alpha\beta}\beta\tns\alpha.
\end{equation*}

The tensor coalgebra structure is given by defining the (reduced)
diagonal $\Delta:T(V)\ra T(V)$ by
\begin{equation}
\Delta(v_1\mtns v_n)=\sum_{k=1}^{n-1}(v_1\mtns
v_k)\tns(v_{k+1}\mtns v_n).
\end{equation}
(We use here the reduced diagonal, because we are not including
the zero degree term in the tensor coalgebra. As a consequence,
$\Delta$ is not injective; in fact, its kernel is $V$.) The tensor
coalgebra is not graded cocommutative under either the \zt\ or the
$\zt\times\Z$-grading, but both gradings are compatible with the
coalgebra structure.

The symmetric coalgebra structure on $\SV$ is given by defining
\begin{equation}
\Delta(v_1\cdots v_n)= \sum_{k=1}^{n-1} \sum_{\sigma\in\sh k{n-k}}
\epsilon(\sigma) v_{\sigma(1)}\cdots v_{\sigma(k)} \tns
v_{\sigma(k+1)}\cdots v_{\sigma(n)}.
\end{equation}
where $\sh pq$ is the set of all \emph{unshuffles} of type
$(p,q)$, that is, the permutations of $p+q$ such that
$\sigma(k)<\sigma(k+1)$ if $k\ne p$.

With this coalgebra structure, and the \zt-grading, $\SV$ is a
cocommutative, coassociative coalgebra.

Similarly, we define the exterior coalgebra structure on
$\bigwedge V$ by
\begin{equation}
\Delta(v_1\mwedge v_n)= \sum\Sb{
k=1\dots n-1\\
\sigma\in\sh k{n-k} } \!\!\!\!\!\! \s\sigma \epsilon(\sigma)
v_{\sigma(1)}\mwedge v_{\sigma(k)} \tns v_{\sigma(k+1)}\mwedge
v_{\sigma(n)}.
\end{equation}
The coalgebra structure on $\bigwedge V$ is coassociative, and is
cocommutative with respect to the $\ztz$-grading.
\section{Coderivations}\label{sect 5}
A coderivation on a graded coalgebra $C$ is a map $d:C\ra C$ such
that
\begin{equation}
\Delta\circ d=(d\tns 1+1\tns d)\circ\Delta.
\end{equation}
Note that this definition depends implicitly on the grading group,
because
\begin{equation*}
(1\tns d)(\alpha\tns\beta)=\s{\ip\alpha d}\alpha\tns d(\beta).
\end{equation*}
The \k-module $\coder(C)$ of graded coderivations has a natural
structure of a graded Lie algebra, with the bracket given by
\begin{equation}
\sb{m,n}=m\circ n-\s{\ip mn}n\circ m.
\end{equation}
A codifferential on a coalgebra $C$ is a coderivation $d$ such
that $d\circ d=0$. We examine the coderivation structure of the
tensor, symmetric and exterior coalgebras.
\subsection{Coderivations of the Tensor Coalgebra}
Suppose that we wish to extend $d_k:V^k\ra V$ to a coderivation
$\hd_k$ of $T(V)$. We are interested in extensions satisfying the
property that $\hd_k(v_1\mcom v_n)=0$ for $n<k$. How this
extension is made depends on whether we consider the \zt\ or the
$\zt\times\Z$-grading. First we consider the \zt-grading, so that
only the parity of $d$ is relevant. Then if we define
\begin{equation*}
\hd_k(v_1\mcom v_n)= \sum_{i=0}^{n-k} \s{(v_1\mplus v_i)d_k}
v_1\mcom v_i, d_k(v_{i+1}\mcom v_{i+k}), v_{i+k+1}\mcom v_n,
\end{equation*}
one can show that $\hd_k$ is a coderivation on $T(V)$ with respect
to the \zt-grading. More generally, one can show that any
coderivation $\hd$ on $T(V)$ is completely determined by the
induced mappings $d_k:V^k\ra V$, and in fact, one obtains
\begin{equation}\label{cder1}
\hd(v_1\mcom v_n)=\bs
\sum \Sb{
1\le k\le n\\
\\
0\le i\le n-k } \bs \s{(v_1\mplus v_i)d_k} v_1\mcom v_i,
d_k(v_{i+1}\mcom v_{i+k}), v_{i+k+1}\mcom v_n.
\end{equation}
Also, one can show that $\hd$ is a codifferential with respect to
the \zt-grading precisely when for $n\ge 1$
\begin{equation*}
\sum \Sb{
k+l=n+1\\
\\
0\le i\le n-k } \s{(v_1\mplus v_i)d_k} d_l(v_1\mcom
v_i,d_k(v_{i+1}\mcom v_{i+k}),v_{i+k+1}\mcom v_n)=0.
\end{equation*}
The module $\coder(T(V))$ of coderivations of $T(V)$ with respect
to the \zt-grading is naturally isomorphic to $\hom(T(V),V)$, so
$\hom(T(V),V)$ inherits a natural structure of a \zt-graded Lie
algebra.

Let us examine the bracket structure on $\hom(T(V),V)$ more
closely. Suppose that for an arbitrary element $d\in\hom(T(V),V)$,
we denote the restriction of $d$ to $V^k$ by $d_k$, and the
extensions of $d$ and $d_k$ as coderivations of $T(V)$ by $\hd$
and $\hd_k$, resp. Also denote by $d_{kl}$ the restriction of of
$\hd_k$ to $V^{k+l-1}$, so that $d_{kl}\in\hom(V^{k+l-1},V^l)$.
The precise expression for $d_{kl}$ is given by equation
(\ref{cder1}) with $n=k+l-1$. It is easy to see that the bracket
of $d_k$ and $\delta_l$ is given by
\begin{equation}
[d_k,\delta_l]= d_k\circ\delta_{lk}-\s{d_k\delta_l}\delta_l\circ
d_{kl}.
\end{equation}
Furthermore, we have $[d,\delta]_n=\sum_{k+l=n+1}[d_k,\delta_l].$
The point here is that $d_{kl}$ and $\delta_{kl}$ are determined
in a simple manner by $d_k$ and $\delta_k$, so we have given a
description of the bracket on $\hom(T(V),V)$ in a direct fashion.
The fact that the bracket so defined has the appropriate
properties follows from the fact that if $\rho=[d,\delta]$, then
$\hat\rho=[\hd,\hdl]$.

Now we consider how to extend a mapping $m_k:V^k\ra V$ to a
coderivation $\hm_k$ with respect to the \ztz-grading. In this
case, note that the bidegree of $m_k$ is given by
$\bd{m_k}=(\e{m_k},k-1)$. The formula for the extension is the
same as before, but with the bidegree in place of the parity. In
other words,
\begin{multline}
\hm_k(v_1\mcom v_n)=\\
\sum_{i=0}^{n-k} \s{(v_1\mplus v_i)m_k+i(k-1)} v_1\mcom v_i,
m_k(v_{i+1}\mcom v_{i+k}), v_{i+k+1}\mcom v_n.
\end{multline}
Similarly, if we consider an arbitrary coderivation $\hm$ on
$T(V)$, then it is again determined by the induced mappings
$m_k:V^k\ra V$, and we see that
\begin{multline}
\hm(v_1\mcom v_n)=\\
\sum \Sb{
1\le k\le n\\
\\
0\le i\le n-k } \bs\s{(v_1\mplus v_i)m_k +i(k-1)} v_1\mcom v_i,
m_k(v_{i+1}\mcom v_{i+k}), v_{i+k+1}\mcom v_n.
\end{multline}
Also, one obtains that $\hm$ is a codifferential with respect to
the $\ztz$-grading is equivalent to the condition that for $n\ge
1$,
\begin{equation}\label{ztzcodiff}
\bs\sum \Sb{
k+l=n+1\\
\\
0\le i\le n-k }\bs \s{(v_1\mplus v_i)m_k +i(k-1)} m_l(v_1\mcom
v_i,m_k(v_{i+1}\mcom v_{i+k}),v_{i+k+1}\mcom v_n)=0.
\end{equation}

The module $\coder(T(V))$ of coderivations of $T(V)$ with respect
to the \ztz-grading is naturally isomorphic to
$\bigoplus_{k=1}^\infty\hom(V^k,V)$, rather than $\hom(T(V),V)$,
because the latter module is the direct product of the modules
$C^k(V)=\hom(V^k,V)$. However, we would like to consider elements
of the form $\hm=\sum_{k=1}^\infty\hm_k$, where $\hm_k$ has
bidegree $(\e{m_k},k-1)$. Such an infinite sum is a well defined
element of $\hom(T(V),T(V))$, so by abuse of notation, we will
define $\coder(T(V))$ to be the module of such infinite sums of
coderivations.
 With this convention, we now have a natural
isomorphism between $\coder(T(V))$ and $\hom(T(V),V)$.
Furthermore, the bracket of coderivations is still well defined,
and we consider $\coder(T(V)$ to be a \ztz-graded Lie algebra. The
reason that the bracket is well defined is that any homogeneous
coderivation has bidegree $(m,n)$ for some $n\ge 0$, so the
grading is given by $\zt\times\N$ rather than the full group \ztz.
In structures where a $\Z$-grading reduces to an \N-grading, it is
often advantageous to replace direct sums with direct products.

Using the same notation convention as in the \zt-graded case, we
note that if $m,\mu\in\hom(T(V),V)$, then we have
\begin{equation}
[m_k,\mu_l]=m_k\circ\mu_{lk}-\s{\ip{m_k}{\mu_l}}\mu_l\circ m_{kl},
\end{equation}
and $[m,\mu]_n=\sum_{k+l=n+1}[m_k,\mu_l]$.

\subsection{Coderivations of the Symmetric Coalgebra}
Suppose that we want to extend $m_k:S^k(V)\ra V$ to a coderivation
$\hm_k$ of $\SV$ such that $\hm_k(v_1\cdots v_n)=0$ for $k<n$.
Define
\begin{equation}
\hm_k(v_1\cdots v_n)= \sum_{\sigma\in\sh k{n-k}} \epsilon(\sigma)
m_k(v_{\sigma(1)}\mcom v_{\sigma(k)})v_{\sigma(k+1)} \cdots
v_{\sigma(n)}.
\end{equation}
Then $\hm_k$ is a coderivation with respect to the \zt-grading. In
general, suppose that $\hm$ is a coderivation on the  symmetric
coalgebra. It is not difficult to see that if $m_k:S^k(V)\ra V$ is
the induced map, then $\hm$ can be recovered from these maps by
the relations
\begin{equation}
\hm(v_1\cdots v_n)= \sum \Sb{
1\le k\le n\\
\\
\sigma\in\sh k{n-k} } \epsilon(\sigma) m_k(v_{\sigma(1)}\mcom
v_{\sigma(k)})v_{\sigma(k+1)} \cdots v_{\sigma(n)}.
\end{equation}
{}From this, we determine that there is a natural isomorphism
between $\coder(\SV)$, the module of coderivations of $V$, and
$\hom(\SV,V)$. Thus $\hom(\SV,V)$ inherits the structure of a
graded Lie algebra. Also, $\hm$ is a codifferential when for all
$n$,
\begin{equation}
\sum \Sb{ k+l=n+1
\\
\sigma\in\sh k{n-k} } \epsilon(\sigma) m_l(m_k(v_{\sigma(1)}\mcom
v_{\sigma(k)}),v_{\sigma(k+1)} \mcom v_{\sigma(n)})=0.
\end{equation}

In general, it is not possible to extend a map $m_k:S^k(V)\ra V$
as a coderivation with respect to the \ztz-grading; the signs
won't work out correctly. Similarly, the \zt-grading is necessary
to have a good theory of coderivations of the symmetric coalgebra.
\subsection{Coderivations of the Exterior Coalgebra}
Suppose that we want to extend $l_k:\bigwedge^k(V)\ra V$ to a
coderivation $\hl_k$ of $\EV$ such that $\hl_k(v_1\mcom v_n)=0$
for $k<n$. Define
\begin{equation}
\hl_k(v_1\mcom v_n)= \sum_{\sigma\in\sh k{n-k}}
\s{\sigma}\epsilon(\sigma) l_k(v_{\sigma(1)}\mcom
v_{\sigma(k)})\wedge v_{\sigma(k+1)} \mwedge v_{\sigma(n)}.
\end{equation}
Then $\hl_k$ is a coderivation with respect to the \ztz-grading.
In general, suppose that $\hl$ is a coderivation on the  exterior
coalgebra. If $l_k:\bigwedge^k(V)\ra V$ is the induced map, then
$\hl$ can be recovered from these maps by the relations
\begin{equation}
\hl(v_1\cdots v_n)= \sum \Sb{
1\le k\le n\\
\\
\sigma\in\sh k{n-k} } \s{\sigma}\epsilon(\sigma)
l_k(v_{\sigma(1)}\mcom v_{\sigma(k)})\wedge v_{\sigma(k+1)}
\mwedge v_{\sigma(n)}.
\end{equation}
{}From this, we determine that there is a natural isomorphism
between $\coder(\EV)$, the module of coderivations of $V$, and
$\hom(\EV,V)$. Thus $\hom(\EV,V)$ inherits the structure of a
graded Lie algebra. Also, $\hl$ is a codifferential when for all
$n$,
\begin{equation}
\sum \Sb{ k+l=n+1
\\
\sigma\in\sh k{n-k} } \s{\sigma}\epsilon(\sigma)
l_l(l_k(v_{\sigma(1)}\mcom v_{\sigma(k)}),v_{\sigma(k+1)} \mcom
v_{\sigma(n)})=0.
\end{equation}
As the \zt-grading is necessary for coderivations of the symmetric
algebra, so the \ztz-grading is essential for a good theory of
coderivations of the exterior coalgebra.

\section{Cohomology of \ainf\ algebras}\label{ainf-cohom}
In \cite{ps2}, a generalization of an associative algebra, called
a strongly homotopy associative algebra, or \ainf\ algebra was
discussed, and cohomology and cyclic cohomology of this structure
was defined. \ainf\ algebras were introduced by J. Stasheff in
\cite{sta1,sta2}. An \ainf\ algebra structure is essentially an
odd codifferential on the tensor coalgebra; an associative algebra
is simply an odd codifferential determined by a single map
$m_2:V^2\ra V$.

If $V$ is a \zt-graded \k-module, then the parity reversion $\Pi
V$ is the same module, but with the parity of elements reversed.
In other words, $(\Pi V)_0=V_1$ and $(\Pi V)_1=V_0$, where $V_0$
and $V_1$ are the submodules of even and odd elements of $V$,
resp. The map $\pi:V\ra \Pi V$, which is the identity as a map of
sets, is odd. There is a natural isomorphism $\eta:T(V)\ra T(\Pi
V)$ given by
\begin{equation}
\eta(v_1\mtns v_n)= \s{(n-1)v_1+(n-2)v_2\mplus v_{n-1}}\pi
v_1\mtns \pi v_n
\end{equation}
Denote the restriction of $\eta$ to $V^k$ by $\eta_k$. Note that
$\eta_k$ is odd when $k$ is odd and even when $k$ is even, so that
$\eta$ is neither an odd nor an even mapping.

Let $W=\Pi V$. Define a bijection between $C(W)=\hom(T(W),W)$ and
$C(V)=\hom(T(V),V)$ by setting
$\mu=\eta^{-1}\circ\delta\circ\eta$, for $\delta\in\CW$. Then
$\mu_k=\eta_1^{-1}\circ\delta_k\circ\eta_k$ and
$\e{\mu_k}=\e{\delta_k}+(k-1)$. In particular, note that if
$\delta_k$ is odd in the \zt-grading, then $\bid(\mu_k)=(k,k-1)$,
so that $\mu_k$ is odd in the \ztz-grading. Now extend
$\delta_k:W^k\ra W$ to a coderivation $\hdl_k$ on $T(W)$ with
respect to the \zt\ grading, so that
\begin{multline}
\hdl_k(w_1\mcom w_n)=\\
\sum_{i=0}^{n-k} \s{(w_1\mplus w_i)\delta_k} w_1\mtns w_i \tns
\delta_k(w_{i+1}\mcom w_{i+k}) \tns w_{i+k+1}\mtns w_n.
\end{multline}
Let $\bmu_k:T(V)\ra T(V)$ be given by $\bmu_k=\eta^{-1}\circ
\hdl_k\circ \eta$. Let $\hmu$ be the extension of $\mu$ as a
\ztz-graded coderivation of $T(V)$. We wish to investigate the
relationship between $\bmu_k$ and $\hmu_k$. For simplicity, write
$w_i=\pi v_i$. Note that $\eta_1=\pi$ is the parity reversion
operator.  So we have $\delta_k=\pi\circ\mu_k\circ\eta_k\inv$.
Thus
\begin{multline}
\bmu_k(v_1\mcom v_n)= \s{r}
\eta^{-1}\delta(w_1\mcom w_n)=\\
\eta^{-1}( \sum_{i=0}^{n-k} \s{r+s_i} w_1\mtns w_i\tns
\delta_k(w_{i+1}\mcom w_{i+k})\tns w_{i+k+1}\mtns w_{n} )
=\\
\eta^{-1}( \sum_{i=0}^{n-k} \s{r+s_i +t_i} w_1\mtns w_i\tns \pi
\mu_k(v_{i+1}\mcom v_{i+k})\tns w_{i+k+1}\mtns w_{n} )
=\\
\sum_{i=0}^{n-k} \s{r+s_i +t_i+u_i} v_1\mtns v_i\tns
\mu_k(v_{i+1}\mcom v_{i+k})\tns v_{i+k+1}\mtns v_{n},
\end{multline}
where
\begin{eqnarray*}
r&=&(n-1)v_1\mplus v_{n-1}\\
s_i&=&(w_1\mplus w_i)\delta_k\\
&=&(v_1\mplus v_i)\mu_k+i(\mu_k+1-k)
+(1-k)(v_1\mplus v_i)\\
t_i&=&(k-1)v_{i+1}\mplus v_{i+k-1}\\
u_i&=&(n-k)v_1\mplus (n-k-i+1)v_i\\
&&+(n-k-i)(\mu_k+v_{i+1}\mplus v_{i+k}) + (n-k-i-1)v_{i+k+1}\mplus
v_{n-1}
\end{eqnarray*}
Combining these coefficients we find that
\begin{equation}
r+s_i+t_i+u_i=(v_1\mplus v_i)\mu_k+i(k-1) +(n-k)\mu_k,
\end{equation}
so that
\begin{multline}\label{hass2}
\bmu_k(v_1\mcom v_n)= \sum_{i=0}^{n-k} \s{ (v_1\mplus
v_i)\mu_k+i(k-1) +(n-k)\mu_k}\\\times v_1\mtns v_i\tns
\mu_k(v_{i+1}\mcom v_{i+k})\tns v_{i+k+1}\mtns v_{n}.
\end{multline}
Thus we see that
\begin{equation}\label{munot}
\bmu_k(v_1\mcom v_n)=\s{(n-k)\mu_k}\hmu_k(v_1\mcom v_n).
\end{equation}
Denote the restriction of $\hmu_k$ to $V^{k+l-1}$ by $\hmu_{kl}$.
Set $n=k+l-1$. Denote
$\bmu_{kl}=\eta_l^{-1}\circ\delta_{kl}\circ\eta_{n}$, so that
$\bmu_{kl}$ is the restriction of $\bmu_k$ to $V^{n}$. Then we
can express equation (\ref{munot}) in the form
$\bmu_{kl}=\s{(l-1)\mu_k}\mu_{kl}$.

More generally, if $\hdl$ is an arbitrary derivation on $T(W)$,
induced by the maps $\delta_k:V^k\ra V$, then it determines maps
$\mu_k:V^k\ra V$ , and $\bmu:T(V)\ra T(V)$, in a similar manner,
and we have
\begin{multline}
\bmu(v_1\mcom v_n)= \sum \Sb{
1\le k\le n\\
\\
0\le i\le n-k } \s{ (v_1\mplus v_i)\mu_k+i(k-1)
+(n-k)\mu_k}\\\times v_1\mtns v_i\tns \mu_k(v_{i+1}\mcom
v_{i+k})\tns v_{i+k+1}\mtns v_{n}.
\end{multline}
The condition that $\hdl$ is a codifferential on $T(W)$ can be
expressed in the form
\begin{equation}
\sum_{k+l=n+1}\delta_l\circ\delta_{kl}=0
\end{equation}
for all $n\ge1$. This condition is equivalent to the condition
\begin{equation}
\sum_{k+l=n+1}\s{(l-1)\mu_k}\mu_l\circ\mu_{kl}=
\sum_{k+l=n+1}\mu_l\circ\bmu_{kl} =0.
\end{equation}
We can express this condition in the form
\begin{multline}
\sum \Sb{ k+l=n+1
\\
0\le i\le n-k } \s{ (v_1\mplus v_i)\mu_k+i(k-1)
+(n-k)\mu_k}\\\times \mu_l(v_1\mcom v_i,\mu_k(v_{i+1}\mcom
v_{i+k}),v_{i+k+1}\mcom v_{n})=0.
\end{multline}
When $\hdl$ is an odd codifferential, $\e{\mu_k}=k$,so the sign in
the expression above is simply $(v_1\mplus v_i)k+i(k-1)+nk-k$.
This is the form in which the signs appear in the definition of an
\ainf\ algebra in \cite{ls,lm}.

An \ainf\ algebra structure on $V$ is nothing more than an odd
codifferential on $T(W).$ This can also be expressed in terms of
the bracket on $\CW$. An odd element $\delta\in\CW$ satisfying
$[\delta,\delta]=0$ determines a codifferential $\hdl$ on $T(W)$.
More precisely, the condition $[\delta,\delta]=0$ is equivalent to
the condition $\hdl^2=0$. Thus $\mu\in\CV$ determines an \ainf\
algebra structure on $V$ when $\delta=\eta\circ\mu\circ\eta^{-1}$
satisfies $[\delta,\delta]=0$. This is not the same condition as
$[\mu,\mu]=0$, nor even the condition $\hmu^2=0$, although we
shall have more to say about this later.

If one considers the situation $V=\Pi W$ instead of $W=\Pi V$,
then an odd codifferential $d$ on $T(W)$ would give rise to an a
mapping $m:\TV\ra V$ satisfying the relations
\begin{multline}\label{mysigns}
\sum \Sb{
0\le 1\le n-k\\
\\
k+l=n+1 } \s{(x_1\mplus x_i)k+i(k-1)+n-k}\\\times m_l(x_1\mcom
x_i,m_k(x_{i+1}\mcom x_{i+k}),x_{i+k+1}\mcom x_{n})=0.
\end{multline}
The signs in the expression above agree with the signs in the
definition of a \ainf\ algebra as given in \cite{ps2,kon}.

{}From these observations, we see that the two sign conventions
for a homotopy associative algebra originate because there are two
natural choices for the relationship between the space $W$, which
carries the structure of an odd differential with respect to the
usual \zt-grading, and $V$, which is its dual. One may choose
either $W=\Pi V$, to get the signs in \cite{ls,lm}, or $V=\Pi W$,
to get the signs in \cite{ps2,kon}. (Actually, one can vary the
definition of $\eta$ to obtain both sets of signs from either one
of these models.)

Let us examine the bracket structure on the space $\coder(T(V))$.
In \cite{gers}, M. Gerstenhaber defined a bracket on the space of
cochains of an associative algebra, which we shall call the
Gerstenhaber bracket. When $V$ is concentrated in degree zero, in
other words, in the non \zt-graded case, the Gerstenhaber bracket
is just the bracket of coderivations, with the \Z-grading. Thus
the Gerstenhaber bracket is given by
\begin{equation}
[\ph_k,\psi_l]=\ph_k\psi_l-\s{(k-1)(l-1)}\psi_l\ph_k,
\end{equation}
for $\ph\in C^k(V)$ and $\psi\in C^l(V)$.

One of the main results of \cite{gers} is that the differential
$D$ of cochains in the cohomology of an associative algebra can be
expressed in terms of the bracket. It was shown that
$D(\ph)=[\ph,m]$ where $m$ is the cochain representing the
associative multiplication. This formulation leads to a simple
proof that $D^2=0$, following from the properties of \Z-graded Lie
algebras. The associativity of $m$ is simply the condition
$[m,m]=0$. Let us recall the proof that an odd homogeneous element
$m$ of a graded Lie algebra satisfying $[m,m]=0$ gives rise to a
differential $D$ on the algebra by defining $D(\ph)=[\ph,m]$. In
other words, we need to show that $[[\ph,m],m]=0$. Recall that $m$
is odd when $\ip{\e m}{\e m}=1$. The graded Jacobi bracket gives
\begin{equation}
[[\ph,m]m] =[\ph[m,m]]+\s{\ip mm}[[\ph,m]m]=-[[\ph,m],m],
\end{equation}
which shows the desired result, in characteristic zero, when the
grading is good. Moreover, we point out that the Jacobi identity
also shows that $D([\ph,\psi])=[\ph,D(\psi)]+\s{\psi
D}[D(\ph),\psi]$, so the differential in the cohomology of an
associative algebra acts as a graded derivation of the Lie
algebra, equipping $\CV$ with the structure of a differential
graded Lie algebra.

We wish to generalize the Gerstenhaber bracket to the \ainf\
algebra case, where we are considering a more general
codifferential $m$ on $T(V)$, in such a manner that the bracket
with $m$ yields a differential graded Lie algebra structure on
$\CV$. If we consider the bracket of coderivations, then a problem
arises when the codifferential is not homogeneous. First of all,
$\hm^2=0$ is not equivalent to the condition $[m,m]=0$. Secondly,
if we define $D(\ph)=[\ph,m]$, then we do not obtain in general
that $D^2=0$. To see this, note that $[m,m]=0$ is equivalent to
$\sum_{k+l=n+1}[m_k,m_l]=0$ for all $n\ge 1$. Let
$\ph_p\in\hom(V^p,V)$. Then
\begin{multline}
[[\ph_p,m],m]_{n+p-1}=
\sum_{k+l=n+1}[[\ph_p,m_k],m_l]=\\
\sum_{k+l=n+1}[\ph_p,[m_k,m_l]]+ \s{\ip{m_k}{m_l}}
[[\ph_p,m_l],m_k]=\\
\sum_{k+l=n+1}\s{k+l+1} [[\ph_p,m_l],m_k]=
\s{n}[[\ph_p,m],m]_{n+p-1}
\end{multline}
Thus we only obtain cancellation of terms when $n$ is odd.
However, this is sufficient to show that in the particular case
where $m_k=0$ for all even or all odd $k$, then $D^2=0$. In this
case we also can show that $\hm^2=0$ is equivalent to $[m,m]=0$ as
well. These problems occur because the product grading on \ztz\ is
not good.

Since these problems do not arise if we are considering \zt-graded
codifferentials (or \Z-graded codifferentials, in the \Z-graded
case), it is natural to consider a codifferential on the parity
reversion(or suspension) of $V$.

For the remainder of this section, let us assume for definiteness
that $W=\Pi V$ and that $d\in\CW$ satisfies $\hd^2=0$. Because $W$
is \zt-graded, $d^2=0$ is equivalent to $[d,d]=0$ for $d$ odd, and
moreover, $\CW$ is a differential graded Lie algebra with
differential $D(\ph)=[\ph,d]$.  Let $m_k=\eta_1\inv\circ
d_k\circ\eta_k$ and $\mu_l=\eta_1\inv\circ\delta_l\circ\eta_l$.
Define a new bracket $\brf$ on $\CV$ by
$\br{m_k}{\mu_l}=\eta_1\inv\circ[d_k,\delta_l]\circ\eta_{k+l-1}$.
It follows easily that
\begin{equation}
\br{m_k}{\mu_l}= \s{(k-1)\mu_l}[m_k,\mu_l].
\end{equation}
Of course this new bracket no longer satisfies the graded
antisymmetry or graded Jacobi identity with respect to the inner
product we have been using on \ztz. However, the bracket does
satisfy these properties with respect to the a different inner
product on \ztz. We state this result in the form of a lemma,
whose proof is straightforward.

\begin{lma}\label{lma1}
Let $V$ be equipped with a \ztz-graded Lie bracket $[\cdot,\cdot]$
with respect to the inner product $\ip{(\bar m,n)}{(\bar
m',n')}=\bar m\bar m'+nn'$ on \ztz. Then the bracket $\brf$ on $V$
given by $\br uv=\s{\deg(u)|v|}[u,v]$ defines the structure of a
\ztz-graded Lie algebra on $V$ with respect to the inner product
$\ip{(\bar m,n)}{(\bar m',n')}=(\bar m+n)(\bar m'+n')$ on \ztz.

The bracket $\brtf$ given by $\brt
uv=\s{\deg(u)(|v|+\deg(v))}[u,v]$ also defines the structure of a
\ztz-graded Lie algebra on $V$ with respect to the second inner
product.
\end{lma}

The bracket $\brf$ on $\CV$ essentially coincides with the bracket
of coderivations on $\CW$, so $[d,d]=0$implies that $\br mm=0$.
Define a differential on $\CV$ by $D(\ph)=\br\ph m$. The homology
of this differential is called the cohomology of the \ainf\
algebra.
 We have been considering
the picture $W=\Pi V$. This means we have been considering \ainf\
algebras as defined in \cite{ls,lm}. If we consider instead the
picture $V=\Pi W$, then the bracket induced on $\CV$ by that on
$\CW$ will be given by $\br
{m_k,\mu_l}=\s{(k-1)(\mu_l+l-1)}[m_k,\mu_l]$, which is the second
modified bracket described in the lemma. Thus we have similar
results. We shall call either one of these two brackets the
modified Gerstenhaber bracket.

The Hochschild cohomology of \ainf\ algebras was defined in
\cite{ps2} as the cohomology given by $D(\ph)=\br \ph m$, and it
was shown that this cohomology classifies the infinitesimal
deformations of an \ainf\ algebra. We shall not go into the
details here. It is important to note however, that unlike the
cohomology theory for an associative algebra, where the \N-grading
on the tensor product gives rise to an  \N-grading on the
cohomology, for an \ainf\ algebra there is only one cohomology
group $H(V)$. (However, $H(V)$ does inherit a natural
\zt-grading.) The reason is that the image of $\ph\in C^n(V)$
under the coboundary operator has a part in all $C^k(V)$ with
$k\ge n$. Only in the case of an associative or differential
graded associative algebra do we get a grading of the cohomology.
\subsection{Cyclic Cohomology of  \ainf\ Algebras}
Now let us suppose that $V$ is equipped with a nondegenerate even
graded symmetric inner product $\ipf$. Graded symmetry means that
$\ip vw=\s{vw}\ip wv$. The inner product induces an isomorphism
between $C^k(V)=\hom(V^k,V)$ and $C^k(V,\k)=\hom(V^{k+1},\k)$,
given by $\ph\mapsto\tilde\ph$, where
\begin{equation}
\tilde\ph(v_1\mcom v_{k+1})=\ip{\ph(v_1\mcom v_k)}{v_{k+1}}.
\end{equation}
Non-degeneracy means that the map $\lambda:\V\ra
\V^*=\hom(\V,\k)$, given by $\lambda(v)(w)=\ip vw$ is an
isomorphism. (When \k\ is a field, this is equivalent to the usual
definition of a non degenerate bilinear form.)

An element $\ph\in\hom(V^k,V)$ is said to be cyclic with respect
to the inner product if
\begin{equation}
\ip{\ph(v_1\mcom v_k)}{v_{k+1}}= \s{k+v_1\ph}
\ip{v_1}{\ph(v_2\mcom v_{k+1})}.
\end{equation}
Then $\ph$ is cyclic if and only if $\tilde\ph$ is cyclic in the
sense that
\begin{equation}
\tilde\ph(v_1\mcom v_{k+1})= \s{n + v_{k+1}(v_1\mplus v_k)}
\tilde\ph(v_{k+1},v_1\mcom v_k).
\end{equation}

If $m\in\CV$, then we say that $m$ is cyclic if $m_k$ is cyclic
for all $k$. If $m$ determines an \ainf\ algebra structure on $V$,
then we say that the inner product is invariant if $m$ is cyclic
with respect to the inner product. Denote the submodule of cyclic
elements in $\CV$ by $CC(V)$, and the submodule of cyclic elements
in $C(V,\k)$ by $CC(V,\k)$. The following lemma shows how to
construct cyclic elements from arbitrary elements of $C^n(V,\k)$.
\begin{lma}\label{lma2}
Suppose that $\tf\in C^n(V,\k)$. Define $C(\tf)\in C^n(V,\k)$ by
\begin{equation}
C(f)(v_1\mcom v_{n+1})= \sum_{0\le i\le n} \s{(v_1\mplus
v_i)(v_{i+1}\mplus v_{n+1})+ni} f(v_{i+1}\mcom v_{i})
\end{equation}
Then $C(\tf)$ is cyclic. Furthermore, $C(\tf)=(n+1)\tf$ if $\tf$
is cyclic.
\end{lma}
The lemma above follows from the technical lemma below, which
simplifies computations with cyclic elements.
\begin{lma}\label{lma3}
If $\tilde f\in C^n(V,\k)$, then $\tilde f$ is cyclic iff whenever
$\alpha=v_1\mtns v_i$ and $\beta=v_i\mtns v_{n+1}$,
\begin{equation}
\tilde f(\alpha\tns\beta)=\s{\alpha\beta+in}\tilde
f(\beta\tns\alpha).
\end{equation}

\end{lma}

The following lemma records the fact that $CC(V)$ is a graded Lie
subalgebra of $\CV$, with respect to the bracket of coderivations.

\begin{lma}\label{lma4}
Let $\ph_k\in C^k(V)$ and $\psi_l\in C^l(V)$. If $\ph$ and $\psi$
are cyclic, then so is $[\ph,\psi]$. Moreover, if $n=k+l-1$, then
\begin{multline}\label{cycrel}
\widetilde{[\ph,\psi]}(v_1\mcom v_{n+1})=\\
\sum \Sb{ 0\le i\le n } \s{(v_1\mplus v_i)(v_{i+1}\mplus
v_{n+1})+in} \tilde\ph_k(\psi_l(v_{i+1}\mcom
v_{i+l}),v_{i+l+1}\mcom v_i),
\end{multline}
where in expressions of this type, indices should be interpreted
$\mod{n+1}$.
\end{lma}

Since $\br{\ph}{\psi}=\s{(k-1)\psi_l}[\ph,\psi]$, it follows that
the modified Gerstenhaber bracket of cyclic elements is also
cyclic. Thus finally, we can state the main theorem, which allows
us to define cyclic cohomology of an \ainf\ algebra.

\begin{thm}\label{thm3}
i)Suppose that $V$ is a \zt-graded \k-module with an inner product
$\ipf$. Suppose that $\ph, \psi\in\CV$ are cyclic. Then
$\br{\ph}{\psi}$ is cyclic. Furthermore, the formula below holds.
\begin{multline}
\widetilde{\br{\ph}{\psi}}(v_1\mcom v_{n+1})=\\
\sum \Sb{
k+l=n+1\\
\\
0\le i\le n }\bs \s{(v_1\mplus v_i)(v_{i+1}\mplus
v_{n+1})+in+(k-1)\psi_l} \tilde\ph_k(\psi_l(v_{i+1}\mcom
v_{i+l}),v_{i+l+1}\mcom v_i).
\end{multline}

Thus the inner product induces a structure of a \ztz-graded Lie
algebra in $CC(V,\k)$  given by
$\br{\tilde\ph}{\tilde\psi}=\widetilde{\br{\ph}{\psi}}$.

ii) If $m$ is an \ainf\ structure on $V$, then there is a
differential in $CC(V,\k)$, given by
\begin{multline}
D(\tilde\ph)(v_1\mcom v_{n+1})=\\
\sum \Sb{
k+l=n+1\\
\\
0\le i\le n }\bs \s{(v_1\mplus v_i)(v_{i+1}\mplus v_{n+1})+in
+(k-1)l} \tilde\ph_k(m_l(v_{i+1}\mcom v_{i+l}),v_{i+l+1}\mcom
v_i).
\end{multline}

iii) If the inner product is invariant, then
$D(\tilde\ph)=\br{\tilde\ph}{\tilde m}$. Thus $CC(V,\k)$ inherits
the structure of a differential graded Lie algebra.
\end{thm}
\begin{proof}
The first statement follows from lemma \ref{lma4}, and the third
follows from the first two. The second assertion follows
immediately from the first when the inner product is invariant.
The general case is a routine verification, which we omit.
\end{proof}
Note that the second statement in the theorem holds even in the
absence of an inner product, because the definition of cyclicity
in $CC(V,\k)$ does not depend on the inner product. Thus cyclic
cohomology of an \ainf\ algebra can be defined independently of
the inner product.

We define $HC(V)$ to be the cohomology associated to the
coboundary operator on $CC(V)$. As in the case of Hochschild
cohomology, the cyclic cohomology is \zt-graded, but does not
inherit an \N-grading, except in the associative algebra case.
\section{Cohomology of \linf\ algbras}\label{linf-cohom}

There is a natural isomorphism $\eta$ between $\bigwedge V$ and
$S(\Pi V)$ which is given by
\begin{equation}
\eta(v_1\mwedge v_n)= \s{(n-1)v_1+(n-2)v_{2}\mplus v_{n-1}}
 \pi v_1\cdots \pi v_n,
\end{equation}
Note that $\eta$ is neither even nor odd. The restriction $\eta_k$
of $\eta$ to $\bigwedge^k V$ has parity $k$. Of course, $\eta$
does preserve the exterior degree. For simplicity in the
following, let $W=\Pi V$ and let $w_i=\pi v_i$, and denote
$C(W)=\hom(S(W),W)$, and $C(V)=\hom(\EV,V)$. We will use
notational conventions as in section \ref{ainf-cohom}, so that for
$d\in C(W)$, $d_k$ will denote the restriction of this map to
$\bigwedge^k W$, $d_{lk}$ will denote the restriction of the
associated coderivation $\hd_l$ to $V^{k+l-1}$ etc.

The following lemma will be useful later on.
\begin{lma}\label{lma5}
Suppose that $\sigma$ is a permutation of $n$ elements. Then
\begin{multline}
\s{(n-1)v_1\mplus v_{n-1}}\s{\sigma}\epsilon(\sigma;v_1\mcom v_n)=\\
\s{(n-1)v_{\sigma(1)}\mplus v_{\sigma(n-1)}}
\epsilon(\sigma;w_1\mcom w_n).
\end{multline}
\end{lma}
\begin{proof}
{}From the properties of the graded exterior algebra, we have
\begin{align*}
\eta(v_{\sigma(1)}\mwedge v_{\sigma(n)})&=
\s{\sigma}\epsilon(\sigma;v_1\mcom v_n) \eta(v_1\mcom v_n)
\\&=\s{\sigma}\epsilon(\sigma;v_1\mcom v_n)
\s{(n-1)v_1\mplus v_{n-1}} w_1\cdots w_n.
\end{align*}
On the other hand, by direct substitution, we have
\begin{align*}
\eta(v_{\sigma(1)}\mwedge v_{\sigma(n)})&=
\s{(n-1)v_{\sigma(1)}\mplus v_{\sigma(n-1)}} w_{\sigma(1)}\cdots
w_{\sigma(n)}\\&= \s{(n-1)v_{\sigma(1)}\mplus v_{\sigma(n-1)}}
\epsilon(\sigma;w_1\mcom w_n)w_1\cdots w_n
\end{align*}
Comparing the coefficients of the two expressions yields the
desired result.
\end{proof}
Let $d\in C(W)$ and  define $l_k=\eta_1\inv\circ d_k\circ\eta_k$.
so that
\begin{equation}
d_k(w_{\sigma(1)}\mcom w_{\sigma(k)}= \s{(k-1)v_{\sigma(1)}\mplus
v_{\sigma(k-1)}} \pi m_k(v_{\sigma(1)}\mcom v_{\sigma(k-1)}).
\end{equation}
Let us abbreviate $\epsilon(\sigma;v_1\mcom v_n)$ by
$\epsilon(\sigma;v)$ and $\epsilon(\sigma;w_1\mcom w_n)$ by
$\epsilon(\sigma;w)$. Let $n=k+l-1$. Define $\bar
l_{lk}=\eta_k\inv\circ d_{lk}\circ\eta_{n}$, where $d_{lk}$ is
given by
\begin{equation}
d_{lk}(w_1\mcom w_n)= \sum_{\sigma\in\sh k{n-k}}
\epsilon(\sigma,w) d_l(w_{\sigma(1)}\mcom w_{\sigma(l)})
w_{\sigma(l+1)}\cdots w_{\sigma(n)}.
\end{equation}
We wish to compute $\bar l_{lk}$ in terms of $l_l$. Now
\begin{multline}
l_{lk}(v_1\mcom v_n)=\\
\sum_{\sigma\in\sh l{n-l}} \!\!\! \epsilon(\sigma,w)
\s{(n-1)v_1\mplus v_{n-1}} \eta_k\inv( d_l(w_{\sigma(1)}\mcom
w_{\sigma(l)}) w_{\sigma(l+1)}\modot w_{\sigma(n)})
\\=
\sum_{\sigma\in\sh l{n-l}} \s{\sigma}\epsilon(\sigma,v) \s{r}
\eta_k\inv( d_l(w_{\sigma(1)}\mcom w_{\sigma(l)})
w_{\sigma(l+1)}\cdots w_{\sigma(n)})
\\=
\sum_{\sigma\in\sh l{n-l}} \s{\sigma}\epsilon(\sigma,v) \s{r+s}
\eta_k\inv( l_l(v_{\sigma(1)}\mcom
v_{\sigma(l)})w_{\sigma(l+1)}\cdots w_{\sigma(n)})
\\=
\sum_{\sigma\in\sh l{n-l}} \s{\sigma}\epsilon(\sigma,v) \s{r+s+t}
l_l(v_{\sigma(1)}\mcom v_{\sigma(l)}) v_{\sigma(l+1)}\cdots
v_{\sigma(n)}
\\=
\sum_{\sigma\in\sh l{n-l}} \s{\sigma}\epsilon(\sigma,v)
\s{(k-1)l_l} l_l(v_{\sigma(1)}\mcom v_{\sigma(l)})
v_{\sigma(l+1)}\cdots v_{\sigma(n)}.
\end{multline}

where
\begin{eqnarray}
r&=&(n-1)v_{\sigma(1)}\mplus v_{\sigma(n-1)}\\
s&=&(l-1)v_{\sigma(1)}\mplus v_{\sigma(l-1)}\\
t&=&(k-1)(l_l+v_{\sigma(1)}\mplus v_{\sigma(l-1)}
(k-2)v_{\sigma(l+1)}\mplus v_{\sigma(n-1)}
\end{eqnarray}

Thus we deduce immediately that $\bar l_{lk}=\s{(k-1)l_l}l_{lk}$,
where $l_{lk}$ is the restriction of the coderivation $\hat l_k$
to $V^{k+l-1}$. This formula is identical to the formula we
deduced in section \ref{ainf-cohom} connecting $\bm_{lk}$ to
$\hm_{lk}$.

Suppose that $d$ is an odd codifferential, so that $d^2=0$. This
is equivalent to $\sum_{k+l=n+1} d_k\circ d_{lk}=0$, which is
equivalent to the relations $\sum_{k+l=n+1} \s{(k-1)l}l_k\circ
l_{lk}=0$, since $\e{l_l}=l$. This last relation can be put in the
form
\begin{multline}
\sum \Sb{
k+l=n+1\\
\sigma\in\sh l{n-l} } \s{\sigma}\epsilon(\sigma,v) \s{(k-1)l}
l_k(l_l(v_{\sigma(1)}\mcom v_{\sigma(l)}), v_{\sigma(l+1)}\mcom
v_{\sigma(n)}))=0.
\end{multline}
We say that the maps $l_k$ induce the structure of an \linf\
algebra, or strongly homotopy Lie algebra on $V$. In \cite{ls},
the sign $k(l-1)$ instead of $(k-1)l$ appears in the definition,
but since $k(l-1)-(k-1)l=n+1$, this makes no difference in the
relations.

If we define $[\ph,\psi]$ to be the bracket of coderivations, then
we can define the modified bracket
$\br{\ph}{\psi}=\s{\deg\ph\e\psi}[\ph,\psi]$. The definition of an
\linf\ algebra can be recast in terms of the bracket. In this
language, $l\in\CV$ determines an \linf\ structure on $V$ when
$\br ll=0$. The cohomology of an \linf\ algebra is defined to be
the cohomology on $\CW$ induced by $l$, in other words,
$D(\ph)=[\ph, l]$. This definition makes $\CV$ a differential
graded Lie algebra, with respect to the second inner product on
\ztz. These results are completely parallel to the \ainf\ case.

In \cite{ps2}, the relationship between infinitesimal deformations
of an \ainf\ algebra and the cohomology of the \ainf\ algebra was
explored. The basic result is that the cohomology classifies the
infinitesimal deformations. Since we did not explore this matter
here for \ainf\ algebras, we shall discuss the parallel result for
\linf\ algebras.

An infinitesimal deformation $l_u$ of an \linf\ algebra is given
by taking $l_u=l+u\lambda$, where $u$ is an infinitesimal
parameter whose  parity  chosen so that $(l_u)_k$ has parity $k$.
It follows that $\e{\lambda_k}=\e{u}+k$. Since $u$ must have fixed
parity, this determines the parity of $\lambda_k$.

The situation is more transparent if we switch to the $W$ picture,
so suppose that $l=\eta\inv\circ d\circ \eta$, and
$\lambda=\eta\inv\circ\delta\circ\eta$. Let $d_u=d +u\delta$, and
let us suppose that $d^2=0$, which is equivalent to $l$ giving an
\linf\ structure on $V$. Then $d_u$ is an infinitesimal
deformation of $d$ if $d_u^2=0$. Since $u$ is an infinitesimal
parameter, $u^2=0$, but we also want the parity of $d_u$ to be
odd, so that $\e u=1- \e\delta$. (We assume here that $\delta$ is
homogeneous.) Now $d_u^2=0$ is equivalent to $d^2+u\delta d+
du\delta=0$. Also, $du=\s{ud}ud=-\s{\delta d}ud$, so this
condition reduces to $[\delta, d]=0$; in other words,
infinitesimal deformations are given by cocycles in the cohomology
of the \linf\ algebra.

Trivial deformations are more complicated in the \linf\ case than
for Lie algebras. Consider $d$ as a codifferential of $T(W)$. Two
codifferentials are said to be equivalent if there is an
automorphism of $T(W)$ which takes one of them to the other. In
the Lie algebra case, one only considers automorphisms of $T(W)$
induced by linear isomorphism of $W$ into itself. One can show
that $d_u$ is a trivial infinitesimal deformation precisely when
$\delta$ is a coboundary. Thus the cohomology of $\CW$ classifies
the infinitesimal deformations of the \linf\ algebra. When we
transfer this back to the $V$ picture, note that $m+\s{u}u\lambda$
is the deformed product associated to $d_u$. The condition for
$l_u$ to be an \linf\ algebra still is that $\br \lambda m=0$. We
sumarize these results in the theorem below.
\begin{thm}
Let $l$ be an \linf\ algebra structure on $V$. Then the cohomology
$H(V)$ of $\CV$ classifies the infinitesimal deformations of the
\linf\ algebra.
\end{thm}

Suppose that $l$ is a Lie algebra structure on $V$. Then the Lie
algebra coboundary operator on $V$ coincides up to a sign with the
\linf\ algebra coboundary operator on $V$. This gives a nice
interpretation of the cohomology of a Lie algebra.
\begin{thm}
Let $l$ be an Lie algebra structure on $V$. Then the Lie algebra
cohomology $H(V)$ of $V$ classifies the infinitesimal deformations
of the Lie algebra into an \linf\ algebra.
\end{thm}
\subsection{Cyclic Cohomology of \linf\ Algebras}

Suppose $V$ is equipped with a nondegenerate graded even symmetric
inner product. An element $\ph\in C^n(V)$ if the tensor
$\tilde\ph:V^n\ra\k$, given by
\begin{equation*}
\tilde\ph(v_1\mcom v_{n+1})=\ip{\ph(v_1\mcom v_n)}{v_{n+1}}
\end{equation*}
is (graded) antisymmetric; \ie, $\tilde\ph\in
CC^n(V,\k)=\hom(\bigwedge^{n+1},\k)$. Note that the antisymmetry
condition is equivalent to the cyclicity condition given in the
definition of cyclicity for \ainf\ algebras, given that $\ph$ is
antisymmetric. Since the inner product is non-degenerate, the map
$\ph\mapsto\tilde\ph$ is an even isomorphism between the submodule
$CC^n(\V)$ of $C^n(\V)$ consisting of cyclic elements, and
$CC^{n}(\V,\k)$.

\begin{thm}
i)Suppose that $V$ is a \zt-graded \k-module with an inner product
$\ipf$. Suppose that $\ph, \psi\in\CV$ are cyclic. Then
$\br{\ph}{\psi}$ is cyclic. Furthermore, the formula below holds.
\begin{multline}
\widetilde{\br{\ph}{\psi}}(v_1\mcom v_{n+1})=\\
\sum \Sb{
k+l=n+1\\
\\
\sigma\in\sh(l,k) } \s{\sigma}\epsilon(\sigma)\s{(k-1)\psi_l}
\tilde\ph_k(\psi_l(v_{\sigma(1)}\mcom v_{\sigma(l)}),
v_{\sigma(l+1)}\mcom v_{\sigma(n+1)}),
\end{multline}
Thus the inner product induces a structure of a \ztz-graded Lie
algebra in the module $CC(V)$ consisting of all cyclic elements in
$\CV)$, by defining
$\br{\tilde\ph}{\tilde\psi}=\widetilde{\br{\ph}{\psi}}$.

ii) If $l$ is an \linf\ structure on $V$, then there is a
differential in $CC(V)$, given by
\begin{multline}
D(\tilde\ph)(v_1\mcom v_{n+1})=\\
\sum \Sb{
k+l=n+1\\
\\
\sigma\in\sh(l,k) } \s{\sigma}\epsilon(\sigma)\s{(k-1)l}
\tilde\ph_k(l_l(v_{\sigma(1)}\mcom v_{\sigma(l)}),
v_{\sigma(l+1)}\mcom v_{\sigma(n+1)}),
\end{multline}

iii) If the inner product is invariant, then
$D(\tilde\ph)=\br{\tilde \ph}{\tilde l}$. Thus $CC(V)$ inherits
the structure of a differential graded Lie algebra.
\end{thm}

We denote the cohomology given by the cyclic coboundary operator
on $CC(V)$ as $HC(V)$. It is interesting to note that
$HC^n(V)\equiv H^{n+1}(V,\k)$, where $H(V,\k)$ is interpreted as
the cohomology of $V$ with coefficients in the trivial module
$\k$. Suppose that $V$ is an \linf\ algebra with an invariant
inner product. Then an infinitesimal deformation $l_t=l+t\ph$
preserves the inner product, that is the inner product remains
invariant under $l_t$, precisely when $\ph$ is cyclic. Thus we see
that cyclic cocycles correspond to infinitesimal deformations of
the \linf\ structure which preserve the inner product. In a
similar manner as before, cyclic coboundaries correspond to
trivial deformations preserving the inner product. Thus we have
the following classification theorem.
\begin{thm}
Let $l$ be an \linf\ algebra structure on $V$, with an invariant
inner product. Then the cyclic cohomology $HC(V)$ classifies the
infinitesimal deformations of the \linf\ algebra preserving the
inner product.
\end{thm}
When $l$ determines a Lie algebra structure on $V$, the cohomology
$HC(V)$ has a natural \Z-grading, and the group $HC^2(V)$ controls
deformations of the Lie algebra structure.
\begin{thm}
Let $l$ be a Lie algebra structure on $V$, with an invariant inner
product. Then the second cyclic cohomology group $HC^2(V)$
classifies the infinitesimal deformations of the Lie algebra
preserving the inner product.
\end{thm}

For a Lie algebra, we have $HC^n(V)\equiv H^{n+1}(V,\k)$, the Lie
algebra cohomology of $V$ with trivial coefficients. Thus the
deformations preserving an invariant inner product are classified
by the cohomology of the Lie algebra with trivial coefficients.
Finally, we give a nice interpretation of the cyclic cohomology of
a Lie algebra.

\begin{thm}
Let $l$ be an Lie algebra structure on $V$, with an invariant
inner product. Then the Lie algebra cyclic cohomology $H(V)$ of
$V$ classifies the infinitesimal deformations of the Lie algebra
into an \linf\ algebra preserving the inner product.
\end{thm}

Note that as in the case of \ainf\ algebras, $C(V)=\hom(\bigwedge
V,V)$ is really the direct product of its graded subspaces
$C^k(V)$, so it is not a graded Lie algebra in the strict sense of
the definition, because it is not the direct sum of the graded
subspaces. Moreover, the cohomology $H(V)$ does not inherit a
natural \Z-grading in general. Similarly, the space $CC(V)$ of
cyclic cochains is a direct product, and $HC(V)$ does not have a
natural \Z-grading. Nevertheless, the good \ztz-grading does equip
these space with a natural \zt-grading, so that the cohomology is
\zt-graded. These issues are addressed in more detail in
\cite{fiapen1}.

It is interesting to note that the cyclic cohomology of \linf\
algebras has more symmetry than that of \ainf\ algebras; the full
symmetric group acts on the space of cyclic cochains for \linf\
algebras. This difference has some consequences in terms of the
actions of these algebras on graph complexes.  In \cite{pen3}, it
is shown that \linf\ algebras with an invariant inner product act
on the ordinary graph complex, while \ainf\ algebras with an
invariant inner product act on the space of ribbon graphs, which
are given by graphs equipped with a cyclic order at each vertex.
\section{Acknowledgements}
The author would like to thank Albert Schwarz, Dmitry Fuchs and
James Stasheff for reading the original version \cite{pen2} of this article
and providing useful suggestions.
%\bibliography{global}

\begin{thebibliography}{10}

\bibitem{aksz}
M.~Alexandrov, M.~Kontsevich, A.~Schwarz, and O.~Zaboronsky,
\emph{The geometry
  of the master equation and topological quantum field theory}, Internat.
  Journ. Modern Phys. \textbf{A12} (1997), 1405--1423.

\bibitem{conn}
A.~Connes, \emph{Non-commutative differential geometry}, Institute
Des Hautes
  Etudes Scientifiques \textbf{62} (1985), 41--93.

\bibitem{fiapen1}
A.~Fialowski and M.~Penkava, \emph{Deformation theory of infinity
algebras},
  preprint math.RT/0101097, 2000.

\bibitem{gers}
M.~Gerstenhaber, \emph{The cohomology structure of an associative
ring}, Annals
  of Mathematics \textbf{78} (1963), 267--288.

\bibitem{getz}
E.~Getzler and J.D.S. Jones, \emph{\hbox{$A_\infty$}-algebras and
the cyclic
  bar complex}, Illinois Journal of Mathematics \textbf{34} (1990), no.~2,
  256--283.

\bibitem{getz2}
\bysame, \emph{Operads, homotopy algebra, and iterated integrals
for double
  loop spaces}, Preprint, 1995.

\bibitem{hoch}
G.~Hochschild, \emph{On the cohomology groups of an associative
algebra},
  Annals of Mathematics \textbf{46} (1945), 58--67.

\bibitem{kast}
D.~Kastler, \emph{Cyclic cohomology within the differential
envelope}, Travaux
  En Cours, 1988.

\bibitem{kon}
M.~Kontsevich, \emph{Feynman diagrams and low dimensional
topology}, First
  European Congress of Mathematics, Paris, 1992, Birkhauser, Basel, 1994,
  pp.~97--121.

\bibitem{lm}
T.~Lada and M.~Markl, \emph{Strongly homotopy {L}ie algebras},
Comm. in Algebra
  \textbf{23} (1995), 2147--2161.

\bibitem{ls}
T.~Lada and J.~Stasheff, \emph{Introduction to sh {L}ie algebras
for
  physicists}, Intern. J. Theor. Phys \textbf{32} (1993), 1087--1103, Preprint
  hep-th 9209099.

\bibitem{lod}
J.~Loday, \emph{Cyclic homology}, Springer-Verlag, 1992.

\bibitem{mar}
M.~Markl, \emph{A cohomology theory for {A}($m$)-algebras and
applications},
  Journal of Pure and Applied Algebra \textbf{83} (1992), no.~6, 141--175.

\bibitem{pen2}
M.~Penkava, \emph{\linf\ algebras and their cohomology},
  Preprint q-alg/9512014.

\bibitem{pen3}
M.~Penkava, \emph{Infinity algebras and the homology of graph
complexes},
  Preprint q-alg/9601018.

\bibitem{ps2}
M.~Penkava and A.~Schwarz, \emph{\hbox{$A_\infty$} algebras and
the cohomology
  of moduli spaces}, Dynkin Seminar, vol. 169, American Mathematical Society,
  1995, pp.~91--107.

\bibitem{ss}
M.~Schlessinger and J.~Stasheff, \emph{The {L}ie algebra structure
of tangent
  cohomology and deformation theory}, Journal of Pure and Applied Algebra
  \textbf{38} (1985), 313--322.

\bibitem{seib}
P.~Seibt, \emph{Cyclic homology of algebras}, World Scientific,
1987.

\bibitem{sta1}
J.D. Stasheff, \emph{On the homotopy associativity of {H}-spaces
{I}},
  Transactions of the AMS \textbf{108} (1963), 275--292.

\bibitem{sta2}
\bysame, \emph{On the homotopy associativity of {H}-spaces {II}},
Transactions
  of the AMS \textbf{108} (1963), 293--312.

\bibitem{sta4}
\bysame, \emph{The intrinsic bracket on the deformation complex of
an
  associative algebra}, Journal of Pure and Applied Algebra \textbf{89} (1993),
  231--235.

\bibitem{sta3}
\bysame, \emph{Closed string field theory, strong homotopy {L}ie
algebras and
  the operad actions of moduli spaces}, Conf. Proc. Lecture Notes Math. Phys.,
  III, Internat. Press, Cambridge, MA, 1994, Perspectives in mathematical
  physics, pp.~265--288.

\bibitem{umb}
R.~Umble, \emph{The deformation complex for differential graded
hopf algebras},
  Preprint, 1994.

\end{thebibliography}
\providecommand{\bysame}{\leavevmode\hbox
to3em{\hrulefill}\thinspace}

\bibliographystyle{amsplain}
\end{document}